%% file: McKay_Proceedings-4.tex
\documentclass[12pt,oneside,english]{amsart}
\usepackage[T1]{fontenc}
\usepackage[utf8]{inputenc}
\usepackage{geometry}
\usepackage{amstext}
\usepackage{amsthm}
\usepackage{amssymb}
\usepackage{stmaryrd}

\makeatletter
\numberwithin{equation}{section}
\numberwithin{figure}{section}
\theoremstyle{plain}
\newtheorem{thm}{\protect\theoremname}[section]
\theoremstyle{remark}
\newtheorem{rem}[thm]{\protect\remarkname}
\theoremstyle{definition}
\newtheorem{problem}[thm]{\protect\problemname}

\usepackage{hyperref}

\usepackage{prettyref}
\newrefformat{chap}{Chapter \ref{#1}}
\newrefformat{sec}{Section \ref{#1}}
\newrefformat{subsec}{Subsection \ref{#1}}

\newrefformat{fig}{Figure \ref{#1}}
\newrefformat{enu}{\ref{#1}}
\newrefformat{eq}{\ref{#1}}

\newrefformat{thm}{Theorem \ref{#1}}
\newrefformat{lem}{Lemma \ref{#1}}
\newrefformat{cor}{Corollary \ref{#1}}
\newrefformat{rem}{Remark \ref{#1}}
\newrefformat{eq}{Equation (\ref{#1})}

\makeatother

\usepackage{babel}
\providecommand{\problemname}{Problem}
\providecommand{\remarkname}{Remark}
\providecommand{\theoremname}{Theorem}

\begin{document}
\title{Open problems in the wild McKay correspondence and related fields}
\author{Takehiko Yasuda}
\address{Department of Mathematics, Graduate School of Science, Osaka University
Toyonaka, Osaka 560-0043, JAPAN}
\email{takehikoyasuda@math.sci.osaka-u.ac.jp}
\thanks{We would like to thank Editage (www.editage.com) for English language editing. This work was supported by JSPS KAKENHI Grant Number JP18H01112.}
\begin{abstract}
The wild McKay correspondence is a form of McKay correspondence in terms
of stringy invariants that is generalized to arbitrary characteristics.
It gives rise to an interesting connection between the geometry of wild
quotient varieties and arithmetic on extensions of local fields.
 The principal purpose of this article is to collect open problems on the wild McKay correspondence,
as well as those in related fields that the author believes are interesting
or important. It also serves as a survey on the present state
of these fields. 
\end{abstract}

\maketitle
\input{macros.tex}

\global\long\def\ind{\mathrm{ind}}%

\tableofcontents{}

\section{Introduction\label{sec:Introduction}}

This article is a substantial update to the author's paper \cite{yasuda2015openproblems}
written in 2015, which has been posted on his personal webpage. The
principal purpose is to collect open problems on the wild McKay correspondence,
as well as those in related fields that the author believes are interesting
or important. It also serves as a survey on the present state
of these fields. Certain problems in this article appeared earlier than
\cite{yasuda2015openproblems} but are scattered in a number of different
papers. Thus, collecting them here in one place would be meaningful,
in particular, for graduate students or young researchers looking
for problems to work on. The author also wrote similar
survey articles \cite{yasuda2017towardmotivic,yasuda2018thewild}\footnote{The article \cite{yasuda2017towardmotivic} is the second
paper on the wild McKay correspondence 
after the first \cite{yasuda2014thepcyclic} and provides the grand design
of the theory of motivic integration over wild Deligne-Mumford stacks
and its potential applications. The article \cite{yasuda2018thewild}
is a survey paper written in Japanese in 2016. Its English version
will be published in the journal ``Sugaku Expositions.'' } in 2013 and 2016, respectively. Since then, there has been considerable
progress in the field; while some problems have been solved, new problems have
also emerged. Therefore, it is worth producing the present paper
for the sake of an update. 

The \emph{wild McKay correspondence} refers to a generalization
of the McKay correspondence via stringy invariants, which originated
in the works of Batyrev \cite{batyrev1999nonarchimedean} and Denef and
Loeser \cite{denef2002motivic}, to positive and mixed characteristics.
The adjective ``wild'' refers to a situation in which the relevant
finite group has an order divisible by the characteristic of the base
field (or a residue field). As is well known, it is a general phenomenon
that problems become harder due to the wildness. The first attempt in
this direction was \cite{yasuda2014thepcyclic} in which
the author studied the linear actions of the group $\ZZ/p\ZZ$ on affine
spaces in characteristic $p$.
\begin{rem}
This text discusses various works on wild quotient singularities from different
perspectives. First, the modular invariant theory studies the algebraic
nature (such as depth, factoriality, and a set of generators) of the
invariant subring $k[x_{1},\dots,x_{n}]^{G}$ associated with a linear
action of a finite group on a polynomial ring. For a thorough treatment
of the modular invariant theory, the reader is referred to \cite{campbell2011modular}.
A geometric approach to wild quotient singularities, especially in
dimension two, perhaps dates back to Artin's works \cite{artin1975wildlyramified,artin1977coverings},
followed by Peskin's \cite{peskin1983quotientsingularities}.
More recent works include \cite{kir`aly2013groupactions,ito2015wildquotient,lorenzini2020moderately},
to mention a few. 
\end{rem}

The wild McKay correspondence gives rise to the interaction of the geometry
of singularities (in particular, from the perspective of birational
geometry) and arithmetic problems such as counting extensions
of a local field and counting rational points over a finite field.
Therefore, it would shed new light on both birational geometry
and such arithmetic problems. In birational geometry, after the
recent achievement in characteristic zero, there appears to be a trend
toward positive/mixed characteristics. However, generally speaking,
it is harder to treat low characteristics compared with dimension. For
details, the reader is referred to \cite{hacon2019theminimal} and the references
therein. The wild McKay correspondence mainly concerns such
a situation of low characteristics and brings hope that the birational
geometry as we have in characteristic zero will be eventually carried
out in arbitrary characteristics to a certain extent. As for the number
theory, the wild McKay correspondence provides a new way of counting
extensions of a local field. Namely, they had been counted with weights
determined by classical invariants such as discriminants. The wild
McKay correspondence produces numerous new weights of geometric origin
(see Section \ref{sec:The-v-function:-Fr=0000F6hlich=002019s}). 
Expanding this story to global fields is also of great interest and leads to
relating the Batyrev-Manin conjecture and the Malle conjecture (see Section \ref{sec:Over-global-fields}).

The remainder of this paper is organized as follows. Section \ref{sec:Brief-guide}
reviews stringy motives and known results on the wild McKay
correspondence. Open problems are discussed starting from Section
\ref{sec:Log-terminal-singularities}. Generally, the problems are defined first, some of which are very vague, and meaning and further details are provided subsequently. Section \ref{sec:Log-terminal-singularities} discusses
the problem when quotient singularities are log terminal.
In Section \ref{sec:Crepant-resolutions}, problems concerning
crepant resolutions of quotient varieties are discussed, including their existence, construction,
and Euler characteristics. Sections \ref{sec:Rationality} and \ref{sec:Duality}
are concerned with rationality and duality of stringy motives, respectively.
The latter is related to the Poincaré duality and both are related
to the existence of a resolution of singularities. Sections \ref{sec:The-v-function:-Fr=0000F6hlich=002019s}
and \ref{sec:The-moduli-space} respectively discuss problems with the $v$-function
(Fröhlich's module resolvent) and the moduli space of $G$-torsors
over the punctured formal disk, which play important
roles in the wild McKay correspondence. Section \ref{sec:Non-linear-action}
is about non-linear actions, that is, general actions other than linear
actions on affine spaces. Section \ref{sec:Generalization-in-various}
considers the generalization of the McKay correspondence and motivic
integration to various directions. Section \ref{sec:Relation-to-other-wild}
briefly discusses the problems on the relation of the wild McKay correspondence
and other theories in terms of wild ramification. Section \ref{sec:The-derived-wild} discusses
the problem of generalizing the derived McKay correspondence to the wild case. Section \ref{sec:Over-global-fields}
 discusses the McKay correspondence over global fields, in particular,
number fields, which is expected to relate the Batyrev-Manin conjecture
and the Malle conjecture. Finally, Section \ref{sec:quasi-etale} discusses
the behavior of stringy motives under quasi-étale maps and application
to local étale fundamental groups.

\subsection*{Convention}

Throughout the paper, an algebraically closed field $k$ is assumed
unless specified otherwise. This is for simplicity, and most arguments
and problems are valid over any base field. The characteristic
of $k$ is denoted by $p$. When considering the action of a finite group $G$
on an affine space $\AA_{k}^{d}$, it is supposed that the action is linear
and has no pseudo-reflection, unless noted otherwise (see Section
\ref{subsec:The-wild-McKay}).

\section{Brief guide to the wild McKay correspondence and stringy motives\label{sec:Brief-guide}}

\subsection{Wild actions and wild quotients\label{subsec:Wild-actions-and}}

Although the term ``McKay correspondence'' can have multiple
mathematical interpretations today, the common theme shared by them is the relation between
a linear representation of a finite group and the associated quotient
variety. Sometimes, more generally, one considers a smooth variety
with a finite group action instead of a linear representation. A finite
group action is said to be \emph{wild} if the base field has a positive
characteristic $p>0$ and the order of the group (or orders of stabilizer/isotropy
groups) is divisible by $p$. Otherwise, an action is said to be \emph{tame.}\footnote{\emph{In fact, there are a number of definitions of tameness/wildness.
See \cite{kerz2009ondifferent}.}}\emph{ }In particular, every action in characteristic zero is tame.
As these terms suggest, wild actions are often much more difficult
to study than tame ones. Similarly, \emph{wild quotient singularities}, that is,
singularities appearing on the quotient variety associated with a wild
action, do not satisfy many good properties of tame quotient singularities.

\subsection{Stringy motives}

An approach to the McKay correspondence that originated in the works
of Batyrev \cite{batyrev1999nonarchimedean} and Denef-Loeser \cite{denef2002motivic}
uses stringy invariants. Let $X$ be a normal $\QQ$-Gorenstein variety
$X$ over a field $k$, and let $\J_{\infty}X$ denote its arc space,
which parametrizes arcs on $X$, that is, morphisms $D:=\Spec k\tbrats\to X$.
For a positive integer $r$ such that $rK_{X}$ is Cartier, we can
define a function 
\[
F_{X}\colon\J_{\infty}X\to\frac{1}{r}\ZZ_{\ge0}\cup\{\infty\}.
\]
Roughly speaking, this function measures the difference in the two
sheaves $\bigwedge^{\dim X}\Omega_{X/k}$ and the invertible sheaf
$\omega_{X}^{[r]}=\cO(rK_{X})$ along each arc.
The \emph{stringy motive} of $X$ is defined by the motivic
integral,
\[
\M_{\st}(X):=\int_{\J_{\infty}X}\LL^{F_{X}}\,d\mu_{X},
\]
which is an element of some version of the complete Grothendieck ring
of varieties (see \cite{yasuda2020motivic}). As usual, $\LL$ denotes
the Lefschetz motive, that is, the class $\{\AA_{k}^{1}\}$ of affine
line in the Grothendieck ring of varieties. (The class of
a variety $Z$ has been denoted by $\{Z\}$ instead of the more standard $[Z]$; Square brackets are reserved for quotient stacks.) This version of the complete
Grothendieck ring of varieties needs to contain the fractional power
$\LL^{1/r}$ of $\LL$. Note that the integral may diverge; in this case, we set $\M_{\st}(X):=\infty$. For a constructible subset $C$,
we can also define the \emph{stringy motive along $C$} as
\begin{equation}
\M_{\st}(X)_{C}:=\int_{(\J_{\infty}X)_{C}}\LL^{F_{X}}\,d\mu_{X}\label{eq:Mst def}
\end{equation}
by restricting the domain of integral to the space $(\J_{\infty}X)_{C}$
of arcs passing through $C$. 

Suppose that there exists a resolution of singularities, $f\colon Y\to X$,
such that the relative canonical divisor $K_{Y/X}=K_{Y}-f^{*}K_{X}$
has simple normal crossing support. Let us write $K_{Y/X}=\sum_{i\in I}a_{i}E_{i}$
with $E_{i}$ exceptional prime divisors and $a_{i}$ nonzero rational
numbers. For each non-empty subset $J\subset I$, we define 
\[
E_{J}^{\circ}:=\bigcap_{j\in J}E_{j}\setminus\bigcup_{i\in I\setminus J}E_{i}.
\]
For $J=\emptyset$, we set $E_{\emptyset}^{\circ}:=Y\setminus\bigcup_{i\in I}E_{i}$.
Then, we have the formula
\begin{equation}
\M_{\st}(X)_{C}:=\sum_{J\subset I}\{E_{J}^{\circ}\cap f^{-1}(C)\}\prod_{j\in J}\frac{\LL-1}{\LL^{1+a_{i}}-1},\label{eq:Explicit}
\end{equation}
if $a_{i}+1>0$ for every $j$ such that $E_{j}$ meets $f^{-1}(C)$.
Otherwise, $\M_{\st}(X)_{C}=\infty$. In particular, assuming the existence
of such a resolution as above, we have the following equivalences:
\[
\M_{\st}(X)\ne\infty\Leftrightarrow\M_{\st}(X)_{X_{\sing}}\ne\infty\Leftrightarrow X\text{ has only log-terminal singularities.}
\]
Here, $X_{\sing}$ denotes the singular locus of $X$.

A resolution of singularities, $f\colon Y\to X$, is called \emph{crepant} if $K_{Y/X}=0$. For a crepant resolution $f\colon Y\to X$, we have
\[
\M_{\st}(X)_{C}=\{f^{-1}(C)\}.
\]
In particular,
\begin{equation}
\M_{\st}(X)=\{Y\}.\label{eq:crep-st}
\end{equation}

\subsection{The wild McKay correspondence\label{subsec:The-wild-McKay}}

Suppose that a finite group $G$ acts on an affine space $\AA_{k}^{d}$
linearly and faithfully. For simplicity, suppose also that $G$ contains
no pseudo-reflection; $g\in G$ is called a \emph{pseudo-reflection}
if the fixed-point locus $(\AA_{k}^{d})^{g}\subset\AA_{k}^{d}$ has
codimension one. In the tame case, the McKay correspondence in terms
of stringy motive is formulated as follows \cite{batyrev1999nonarchimedean,denef2002motivic}.
\begin{equation}
\M_{\st}(\AA_{k}^{d}/G)=\sum_{[g]\in\Conj(G)}\LL^{d-\age(g)}\label{eq:tame}
\end{equation}
Here, $\Conj(G)$ is the set of conjugacy classes of $G$; $\age(g)$
is the age of $g$, a basic invariant in the McKay correspondence.
In particular, this equality shows that if $Y\to\AA_{k}^{d}/G$ is
a crepant resolution, then we have
\begin{equation}
\chi(Y)=\sharp\Conj(G).\label{eq:EulerConj}
\end{equation}
Here, $\chi(Y)$ is the Euler characteristic of $Y$ defined with either
singular cohomology (the case $k=\CC$) or $l$-adic cohomology (the
general case). Equality \eqref{eq:tame} is generalized to the wild
case as follows \cite{yasuda2020motivic}.
\begin{equation}
\M_{\st}(\AA_{k}^{d}/G)=\int_{\Delta_{G}}\LL^{d-v}\label{eq:wild}
\end{equation}
Here, $\Delta_{G}$ is the moduli space of $G$-torsors over the punctured
formal disk, $\Spec k\tpars$, and $v$ is a locally constructible
function $\Delta_{G}\to\frac{1}{\sharp G}\ZZ_{\ge0}$. Broadly,
the space $\Delta_{G}$ may be infinite-dimensional but is always
the disjoint union of at-most countably many varieties, and the last
integral is defined to be the sum
\[
\sum_{a\in\frac{1}{\sharp G}\ZZ_{\ge0}}\{v^{-1}(a)\}\LL^{d-a}.
\]
For more details regarding the definition of this integral, see \cite{tonini2019moduliof}.
In the tame case, $\Delta_{G}$ has only finitely many points and
is identified with $\Conj(G)$, which allows us to regard \eqref{eq:wild}
as a generalization of \eqref{eq:tame}. We may think of the left-hand side of \eqref{eq:wild} as an invariant concerning the geometry
of $\AA_{k}^{d}/G$. However, the right-hand side has a more
arithmetic nature. The integral is an analog of weighted count of
the Galois extensions of $\QQ_{p}$. 

The same integral is also regarded as the stringy motive $\M_{\st}([\AA_{k}^{d}/G])$
of the quotient stack $[\AA_{k}^{d}/G]$, and equality \eqref{eq:wild}
is rephrased as 
\begin{equation}
\M_{\st}(\AA_{k}^{d}/G)=\M_{\st}([\AA_{k}^{d}/G]).\label{eq:wild-rephrased}
\end{equation}
Note that since $[\AA_{k}^{d}/G]$ is a smooth Deligne-Mumford stack
and the morphism $[\AA_{k}^{d}/G]\to\AA_{k}^{d}/G$ is quasi-finite,
proper, and birational, this morphism is a crepant resolution in the
category of Deligne-Mumford stacks. Equality \eqref{eq:wild-rephrased}
is then viewed as a special case of the more general result that the
stringy motive is invariant under crepant proper maps \cite{yasuda2020motivic}:
if $(\cY,E)\to(\cX,D)$ is a crepant proper birational morphism of
``stacky log pairs,'' then 
\[
\M_{\st}(\cY,E)=\M_{\st}(\cX,D).
\]

\section{Log-terminal singularities\label{sec:Log-terminal-singularities}}
\begin{problem}
\label{prob:log term}When are wild quotient singularities log terminal
(resp.~terminal, canonical, log canonical)?
\end{problem}

The four important properties of singularities in the minimal model program
are terminal, canonical, log terminal, and log canonical. There are
the following implications among them:
\[
\text{terminal}\Rightarrow\text{canonical}\Rightarrow\text{log terminal}\Rightarrow\text{log canonical}
\]
Tame quotient singularities are always log terminal. Wild quotient
singularities are sometimes log terminal, but sometimes they are not. It is desirable
to have useful criteria to decide whether the given wild quotient
singularity is log terminal or not. Focusing on linear actions $G\curvearrowright\AA_{k}^{d}$,
we may ask the following more specific question:
\begin{problem}
\label{prob:rep crit}Find a purely representation-theoretic criterion
for whether $\AA_{k}^{d}/G$ is log terminal (resp.~terminal, canonical,
log canonical). 
\end{problem}

Note that for terminal and canonical singularities, the problem makes
sense in the tame case as well. The criterion in this case is known as
the Reid--Shepherd-Barron--Tai criterion (e.g., see \cite[Cor.\ 6]{yasuda2006motivic}),
which uses ages of $g\in G\setminus\{1\}$, invariants determined
by eigenvalues. In the wild case, we have a criterion for each of
the above four classes of singularities in terms of the integral $\int_{\Delta_{G}}\LL^{d-v}$
in \eqref{eq:wild} (see \cite[Cor\ 1.4]{yasuda2020motivic}). However,
to compute this integral, we have to compute the moduli space $\Delta_{G}$
and the function $v$ on it, which is generally complex and not
representation-theoretic. When $G=\ZZ/p^{n}\ZZ$, we have an affirmative
answer to the problem \cite{yasuda2014thepcyclic,yasuda2019discrepancies,tanno2020thewild,tanno2021onconvergence}. 

Yamamoto \cite{yamamoto2021pathological2} proved that the quotient
variety $\AA_{k}^{3}/(\ZZ/3\ZZ)^{2}$ in characteristic three is not
log canonical, provided that the group has no pseudo-reflection. However,
from any proper subgroup $G\subsetneq(\ZZ/3\ZZ)^{2}$, based on \cite{yasuda2014thepcyclic},
the quotient $\AA_{k}^{3}/G$ is canonical. This example shows that
we cannot reduce the problem to the case of cyclic groups by restricting
the given representation to cyclic subgroups. 

Problem \ref{prob:log term} would be much harder when $G$ non-linearly
acts on a smooth variety and when $G$ acts on a singular variety.
Note that the tame action on a smooth variety is always locally linearizable;
each point $x\in X$ has local coordinates for which the action of
the stabilizer is linear. This is no longer the case for the wild action.
Perhaps we should start with the actions of $\ZZ/p$ on $\Spec k\llbracket x,y\rrbracket$
fixing only the closed point. This action is never linear. Artin \cite{artin1975wildlyramified}
and Peskin \cite{peskin1983quotientsingularities} gave partial results
to Problem \ref{prob:log term} in characteristics two and three, respectively.
Looking at their computation, there appear to be two important numerical
invariants of an action. The first one is the number $r$ of Jordan
blocks for the linear action on $(x,y)/(x,y)^{2}$. The second one
is the largest integer $n$ such that, after a suitable choice of
coordinates $x$ and $y$, the action on $k\llbracket x,y\rrbracket$
is linear modulo $(x,y)^{n+2}$; this is an analog of ramification
jump in the context of local fields. 
\begin{problem}
Does $(\Spec k\llbracket x,y\rrbracket)/G$ being log terminal depend only on the above pair $(r,n)$?
\end{problem}

It might be too naive to expect this to be the case. If that is the case, we may
look for more numerical invariants to get an affirmative result.
\begin{problem}
Even when both sides of \eqref{eq:wild} diverge, we may obtain an equality
of motivic elements (or an equality of more meaningful quantities
instead of the equality $\infty=\infty$) by ``renormalizing'' them.
\end{problem}

If $X$ is not log terminal, then $\M_{\st}(X)=\infty$ by definition.
We can still define the dimension of $\M_{\st}(X)$ in any case; however,
it is also infinite if $X$ is not log canonical. Equality \eqref{eq:wild}
holds even if $\AA_{k}^{d}$ is not log canonical. However, in this case,
it only says that both sides are infinity and have infinite dimension.
Can we extract certain finite values by renormalization? As a step in
this direction, Veys' work \cite{veys2003stringy,veys2004stringy}
on stringy invariants of varieties with non-log-canonical singularities
in characteristic zero is notable. 
\begin{rem}
A possible approach to this problem may be to use analytic continuation.
Suppose that we have a linear action $V=\AA_{k}^{d}$ of a finite
group $G$ such that $\int_{\Delta_{G}}\LL^{d-v}=\infty$. For an
integer $s>0$, the direct sum $V^{\oplus s}=\AA_{k}^{sd}$ of $s$
copies of $G$ as a $G$-representation leads to $\int_{\Delta_{G}}\LL^{d-sv}$.
The larger the value of $s$, the more likely the last integral converges. Now
imagine that $s$ is allowed to vary as a complex number. We hope
that the integral may still make sense and converge for $s$ with
a sufficiently large real part. We then try to take its analytic continuation
to the entire complex $s$-plane. 
\end{rem}

\section{Crepant resolutions\label{sec:Crepant-resolutions}}
\begin{problem}
When does there exist a crepant resolution of $\AA_{k}^{d}/G$? Is
it possible to construct a crepant resolution as a certain moduli
space?
\end{problem}

The McKay correspondence, \eqref{eq:tame} and \eqref{eq:wild}, takes
the simplest form when there exists a crepant resolution $Y\to\AA_{k}^{d}/G$,
owing to \eqref{eq:crep-st}. We would like to know when this is the
case. In characteristic zero, for a finite subgroup $G\subset\mathrm{SL}_{d}(\CC)$
with $d\le3$, there exists a crepant resolution of $\AA_{\CC}^{d}/G$.
In dimension two, the minimal resolution is a crepant resolution.
In dimension three, this was proved via a case-by-case analysis by
Roan, Ito, and Markushevich (see \cite{roan1996minimal} and references
therein). It was later proved in \cite{nakamura2001hilbert,bridgeland2001themckay}
that the $G$-Hilbert scheme becomes a crepant resolution in this
case, thereby giving a moduli-theoretic construction of a crepant resolution. 

For the wild case, the following is a list of all the examples that the author
knows:
\begin{enumerate}
\item The symmetric product of an affine plane $\AA_{k}^{2n}/S_{n}=S^{n}\AA_{k}^{2}$
in any characteristic; the Hilbert scheme of $n$ points on $\AA_{k}^{2}$,
$\mathrm{Hilb}^{n}(\AA_{k}^{2})$, gives a crepant resolution \cite{kumar2001frobenius,brion2005frobenius}.
\item A quotient of $\AA_{k}^{2n}$ by the wreath product $(\ZZ/2\ZZ)\wr S_{n}$
in characteristic $\ne2$, as in \cite{wood2017massformulas}.
\item Quotients of $\AA_{k}^{3}$ in characteristic three by groups $\ZZ/3\ZZ$
\cite{yasuda2014thepcyclic} and by more general groups, including
$S_{3},A_{4,}S_{4}$ \cite{yamamoto2018acounterexample,yamamoto2021pathological2}.
\item Quotients of $\AA_{k}^{2}$ in characteristic three by an action of
$\ZZ/6\ZZ$ that has pseudo-reflections \cite{chen2020modular}. 
\end{enumerate}
When $G=\ZZ/p\ZZ$, there is an integer invariant denoted by $D_{V}$
of the given representation $V=\AA_{k}^{d}$ (see \cite{yasuda2014thepcyclic}).
A necessary condition for the existence of a crepant resolution of
$X=\AA_{k}^{d}/G$ is $D_{V}=p$ \cite[Cor.\ 6.21]{yasuda2014thepcyclic}.
It is natural to ask:
\begin{problem}
Is the equality $D_{V}=p$ also a sufficient condition? 
\end{problem}

The simplest example with $D_{V}=p$ is the following: 
\begin{problem}
Suppose that $G=\ZZ/p\ZZ$ acts on $\AA_{k}^{2p}$ by the Jordan matrix
\[
\begin{pmatrix}1 & 1\\
0 & 1
\end{pmatrix}^{\oplus p},
\]
which satisfies $D_{V}=p$. Does the associated quotient $\AA_{k}^{2p}/G$
have a crepant resolution?
\end{problem}

The last quotient variety $\AA_{k}^{2p}/G$ appears (at least as per
the author) as the wild counterpart of the cyclic quotient singularity
in characteristic zero of type $\tfrac{1}{l}(1,\dots,1)$ of dimension
$l$. The latter admits a crepant resolution.

For the cyclic group $G=\ZZ/p^{n}\ZZ$ of order $p^{n}$, Tanno \cite{tanno2021onconvergence}
defined a series of invariants, $D_{V}^{(0)},D_{V}^{(1)},\dots,D_{V}^{(n-1)}$,
generalizing the above invariant $D_{V}$ for the group $\ZZ/p\ZZ$.
He used them to address Problem \ref{prob:rep crit} for this group.
\begin{problem}
Is there a necessary condition expressed by using $D_{V}^{(0)},D_{V}^{(1)},\dots,D_{V}^{(n-1)}$
for the existence of a crepant resolution of $\AA_{k}^{d}/(\ZZ/p^{n}\ZZ)$?
\end{problem}

In the tame case, the Euler characteristic of a crepant resolution
of $\AA_{k}^{d}/G$ is equal to the number of conjugacy classes of
$G$ (see \eqref{eq:EulerConj}). This is not true in the wild case, as
Yamamoto's examples show \cite{yamamoto2018acounterexample, yamamotoCrepant} .
There is also a similar example such that the group has pseudo-reflections
\cite{chen2020modular}. 
\begin{problem}
Find a formula for the Euler characteristic of a crepant resolution
of $\AA_{k}^{d}/G$. 
\end{problem}

Note that from \eqref{eq:crep-st}, the two crepant resolutions of the
same variety have the same Euler characteristic. It would also be interesting to study the following two problems, on which little research
has been carried out for the wild case. 
\begin{problem}[cf.~\cite{yamagishi2015crepant}]
Given a quotient variety $\AA_{k}^{d}/G$, count crepant resolutions
of it. When is it unique?
\end{problem}

\begin{problem}[cf.~\cite{craw2004flopsof}]
Given two different crepant resolutions of the same quotient variety
that are both constructed moduli-theoretically, can we connect them
by wall-crossing? 
\end{problem}

\section{Rationality\label{sec:Rationality}}
\begin{problem}
Suppose that $\M_{\st}(\text{\ensuremath{\AA}}_{k}^{d}/G)\ne\infty$.
Is $\M_{\st}(\AA_{k}^{d}/G)$ a rational function in $\LL^{1/\sharp G}$? 
\end{problem}

In the tame case, from \prettyref{eq:wild}, $\M_{\st}(\AA_{k}^{d}/G)$
is the sum of powers of $\LL^{1/\sharp G}$, in particular, a polynomial
in $\LL^{1/\sharp G}$. In the wild case, it is not necessarily a
polynomial (see computation in \cite{yasuda2014thepcyclic}). However,
to the best of the author's knowledge, it is always a rational function. 
\begin{rem}
Even in the tame case, if $k$ is not algebraically closed, then $\M_{\st}(\AA_{k}^{d}/G)$
is generally not a rational function in $\LL^{1/\sharp G}$ (cf.~\cite[Th.\ 5.9 and Cor.\ 5.11]{wood2015massformulas}).
However, a weaker rationality as in the next problem holds over any field
in the tame case. 
\end{rem}

\begin{problem}
Suppose that $\M_{\st}(\text{\ensuremath{\AA}}_{k}^{d}/G)\ne\infty$.
Can $\M_{\st}(\AA_{k}^{d}/G)$ always be written as a finite sum of
the form 
\begin{equation}
\sum_{i=1}^{n}\frac{\{Z_{i}\}}{\LL^{a_{i}}-1}\label{eq:rat ex}
\end{equation}
with $a_{i}\in\frac{1}{\sharp G}\ZZ_{>0}$ and varieties $Z_{i}$? 
\end{problem}

If $\AA_{k}^{d}/G$ has a log resolution, then formula \eqref{eq:Explicit}
provides an affirmative answer to this question. In other words,
if the last problem has a negative answer for some quotient variety,
then this quotient variety does not admit any log resolution. The
duality discussed in the next section is another example of properties
that hold if some version of resolution of singularities is available. 

\section{Duality\label{sec:Duality}}
\begin{problem}
\label{prob:dual}Let $o\in\AA_{k}^{d}/G$ be the image of the origin
of $\AA_{k}^{d}$ and let $\M_{\st}(\AA_{k}^{d}/G)^{\vee}$ denote
the ``dual'' of $\M_{\st}(\AA_{k}^{d}/G)$. Does the following equality
holds?
\begin{equation}
\M_{\st}(\AA_{k}^{d}/G)_{o}=\M_{\st}(\AA_{k}^{d}/G)^{\vee}\cdot\LL^{d}\label{eq:dual}
\end{equation}
\end{problem}

Specifically, we need to take certain realization of motives for the dual $\M_{\st}(\AA_{k}^{d}/G)^{\vee}$ to make sense. For
example, let us choose here the Poincaré polynomial realization. The
Poincaré polynomial $\P(X)$ of a variety $X$ is an element of the
polynomial ring $\ZZ[T]$. It defines a ring map 
\begin{equation}
\P\colon\Kzero(\Var_{k})\to\ZZ[T],\label{eq:P}
\end{equation}
which, in particular, sends $\LL$ to $T^{2}$. For a smooth proper
variety $X$, we have
\[
\P(X)=\sum_{i=0}^{2\dim X}(-1)^{i}b_{i}T^{i},
\]
where $b_{i}$ is the $i$-th Betti number with respect to $l$-adic
cohomology, and the Poincaré duality is expressed as the functional
equation
\[
\P(X)=\P(X)^{\vee}\cdot T^{2\dim X}.
\]
Here, the dual $\P(X)^{\vee}$ is the Laurent polynomial obtained by
substituting $T$ in $\P(X)$ with $T^{-1}$. For a log-terminal variety
$X$ admitting a log resolution, sending $\M_{\st}(X)$ by the natural
extension of map \eqref{eq:P} gives the ``stringy Poincaré series''
\[
\P_{\st}(X)\in\ZZ\llparenthesis T^{-1/r}\rrparenthesis.
\]
This is a rational function in $T^{1/r}$, and we can define its dual
$\P_{\st}(X)^{\vee}$. Batyrev's argument \cite{batyrev1998stringy}
shows 
\begin{equation}
\P_{\st}(X)=\P_{\st}(X)^{\vee}\cdot T^{2\dim X}.\label{eq:Poincare}
\end{equation}
The duality in Problem \ref{prob:dual} is related to the Poincaré
duality for the quotient space $(X-\{o\})/\GG_{m}$. See \cite{yasuda2014thepcyclic,wood2017massformulas}
for more details. 

Note that, similar to equality \eqref{eq:wild}, the stringy motive $\M_{\st}(\AA_{k}^{d}/G)_{o}$
at $o$ is expressed as an integral $\int_{\Delta_{G}}\LL^{w}$ with
a function $w$ that is similar to $v$. In the tame case, the duality \eqref{eq:dual}
is a direct consequence of the following equality, which easily follows
from the definition of age: 
\[
\age(g^{-1})=\codim((\AA_{k}^{d})^{g},\AA_{k}^{d})-\age(g).
\]
Here, $(\AA_{k}^{d})^{g}$ denotes the fixed-point locus for the action
of $g$. However, in the wild case, the duality, if true, would be related
to deeper natures of the moduli space $\Delta_{G}$ and of functions
$v$ and $w$ on it. 

\section{The $v$-function/Fröhlich's module resolvent\label{sec:The-v-function:-Fr=0000F6hlich=002019s}}

We now explain the background and motivation for the $v$-function. This
function, which appears in the wild McKay correspondence, \eqref{eq:wild},
is important, since it is a common generalization of the two fundamental
invariants, age and the Artin conductor \cite{wood2015massformulas}.
It is also a special case of Fröhlich's module resolvent \cite{frohlich1976moduleconductors}
(see also \cite[Remark 9.4]{tonini2019moduliof}). 

The $v$-function is a function $\Delta_{G}\to\tfrac{1}{\sharp G}\ZZ_{\ge0}$
associated with a representation $G\to\GL_{d}(k\tbrats)$ over $k\tbrats$,
in particular, to a representation over $k$. Here, $\Delta_{G}$ is
the moduli space of $G$-torsors over the punctured formal disk $D^{*}=\Spec k\tpars$
(see Section \ref{sec:The-moduli-space}). The function is defined
by using the \emph{tuning module} associated with the given
representation and a $G$-torsor over $D^{*}$ \cite[Def.\ 3.1]{wood2015massformulas}
(cf.~\cite[Rem.\ 4.1]{yasuda2016wildermckay}). 

The $v$-function can be viewed as a measure of how much the (normal)
$G$-cover of the formal disk $D=\Spec k\tbrats$ corresponding to
a given $G$-torsor over $D^{*}$ is ramified, and one can use it as
a weight when counting $G$-covers of $D$. In fact, the ``w'' in the
$w$-function, an elder sibling of the $v$-function, stands for ``weight''
(see \cite[Def.\ 6.5]{yasuda2017towardmotivic}), and ``v'' was chosen
because it is the alphabet next to ``w.'' See Section \ref{sec:Over-global-fields}
for further details on this viewpoint of counting applied to global fields. Later studies
\cite{wood2015massformulas,wood2017massformulas,yasuda2016wildermckay,yasuda2020motivic}
gradually showed that $v$ is more fundamental than $w$. 
\begin{problem}
Find a ``nice'' formula for the function $v$ in terms of the basic
invariants of the representation $\rho$ and $G$-torsors. 
\end{problem}

As mentioned above, the function $v$ is the same as age in the tame
case and is expressed in terms of the Artin conductor or the discriminant
in the case where $\rho$ is a permutation representation. However, it
is complex to compute this function for a general representation.
When $G=\ZZ/p^{n}\ZZ$, we have a formula for $v$ in terms of ramification
filtration \cite{yasuda2014thepcyclic,tanno2020thewild}, though it
is rather involved. However, even in the case of group $(\ZZ/p\ZZ)^{2}$,
the function is not determined by ramification filtration \cite{yamamoto-v-func}.
We can consider the same problem as above in mixed characteristic.
There has been little research on this case, except the tame case and the case
of permutation representations (see also Problem \ref{prob:motivic-mixed}). 
\begin{problem}
Is the $v$-function upper-semicontinuous?
\end{problem}

To make this problem precise, we first need to decide what geometric
structure to give to the moduli space $\Delta_{G}$. Thus, the problem
is related to the construction of this moduli space, which we discuss
in the next section. 

The $w$-function mentioned above is used in the following variant
of the wild McKay correspondence:
\[
\M_{\st}(\AA_{k}^{d}/G)_{o}=\int_{\Delta_{G}}\LL^{w}
\]
Cf.~Section \ref{sec:Duality}. This version is not explicitly stated
in the literature, but its point-counting version appears in \cite{yasuda2017thewild}.
We can prove the above equality by combining arguments of \cite{yasuda2017thewild}
and the proof of \eqref{eq:wild} in \cite{yasuda2020motivic}. In
the above equality, we need to use the definition of $w$ in \cite{yasuda2017thewild},
which is slightly different from the original one in \cite{yasuda2017towardmotivic,wood2015massformulas}.
The two definitions coincide for a few important cases (see \cite[Rem.\ 8.2 and Lem.\ 8.3]{yasuda2017thewild}).
However, no example in which they are, indeed, different functions has been found so far. 
\begin{problem}
Do the two definitions of $w$ eventually coincide? If not, construct
a counterexample. 
\end{problem}

\section{The moduli space $\Delta_{G}$\label{sec:The-moduli-space}}
\begin{problem}
Can we construct the moduli space $\Delta_{G}$ of $G$-torsors over
the punctured formal disk $D^{*}=\Spec k\tpars$ as an inductive limit
of Deligne-Mumford stacks of finite type?
\end{problem}

Roughly speaking, for a finite group $G$, we consider the moduli
functor 
\begin{align*}
(\text{Affine schemes}/k) & \to(\text{Sets})\\
\Spec R & \mapsto\left\{ G\text{-torsors over \ensuremath{\Spec R\tpars}}\right\} /\cong
\end{align*}
and ask whether or not this is representable by  a
scheme or a generalization of it. Since $G$-torsors have nontrivial automorphisms, as is usual
in the moduli problem, the fine moduli space should be a stack if it
exists. When group $G$ is a semidirect product $H\rtimes C$
of a $p$-group $H$ and a tame cyclic group $C$, it was proved in
\cite{tonini2020moduliof} that the moduli stack $\Delta_{G}$ is
an inductive limit of Deligne-Mumford stacks of finite type. For a
general finite group, we can construct the ``moduli space'' of $\Delta_{G}$
as a rough geometric structure called \emph{P-scheme} \cite{tonini2019moduliof}.
However, the above problem itself is open. 
\begin{problem}
\label{prob:moduli mix}Construct the moduli space $\Delta_{G}$ for
a local field of characteristic zero such as $\QQ_{p}$. 
\end{problem}

Note that the moduli space should be defined over the residue field
rather than the local field itself. In other words, we are interested
in the following moduli functor:

\begin{align*}
(\text{Affine schemes}/k) & \to(\text{Sets})\\
\Spec R & \mapsto\left\{ G\text{-torsors over \ensuremath{\Spec(K\otimes_{W(k)}W(R))}}\right\} /\cong
\end{align*}
Here, $K$ is the given local field, $k$ is its residue field, and
$W(-)$ denotes the ring of the Witt vectors.
\begin{problem}
\label{prob:Delta general G}Construct the moduli space $\Delta_{G}$
for more general group schemes $G$. 
\end{problem}

For example, for $G=\alpha_{p}$, the moduli space should be the ind-pro
limit of finite dimensional spaces (see \cite{tonini2019noteson}).
This would add extra theoretical complexity in applying the
moduli space to the McKay correspondence and motivic integration (see
Problem \ref{prob:alg gp}).

\section{Non-linear actions\label{sec:Non-linear-action}}

In the tame case, every finite group action on a smooth variety is
locally linear, as was mentioned in Section \ref{sec:Log-terminal-singularities}. This is no longer true in the wild case. The wild McKay correspondence
can be generalized to the non-linear case as follows. 

Suppose that a finite group $G$ faithfully acts on a normal variety
$V$ of dimension $d$. Let $A$ be a $\QQ$-divisor on $V/G$ such
that $K_{V/G}+A$ is $\QQ$-Cartier. For each $G$-torsor $E$ over
$\Spec k\tpars$, we can construct the \emph{normalized untwisting
scheme, $V^{|E|,\nu}$}, which is a $k\tbrats$-scheme of relative
dimension $d$ with a natural morphism $V^{|E|,\nu}\to(V/G)\otimes_{k}k\tbrats$
and carries a natural $\QQ$-divisor $B_{E}$. The McKay correspondence
can then be formulated as
\begin{equation}
\M_{\st}(V/G,A)=\int_{\Delta_{G}}\M_{\st}(V^{|E|,\nu},B_{E})/\Aut(E),\label{eq:non-lin}
\end{equation}
where $\Aut(E)$ is the automorphism group of $E$ as a $G$-torsor.
See \cite[Conj.\ 1.3]{yasuda2016wildermckay} and \cite[Th.\ 13.3 and Cor.\ 13.4]{yasuda2020motivic}. 
\begin{problem}
Compute the right-hand side of \eqref{eq:non-lin}. In particular,
compute normalized untwisting schemes $V^{|E|,\nu}$. 
\end{problem}

The only case where these have been performed is the one where $V$ is a smooth
curve in characteristic two, and the group $G$ has order two and acts
on $V$ in the mildest way \cite[Section 10]{yasuda2016wildermckay}.
In this case, there appeared an infinite series of rational double
points on surfaces $V^{|E|,\nu}$ from Artin's classification \cite{artin1977coverings}
as $G$-torsor $E$ varies. 
\begin{problem}
If $V$ is a smooth curve, then what singularities do surfaces $V^{|E|,\nu}$
have? Does every type of rational double point appear on some $V^{|E|,\nu}$?
\end{problem}

\begin{problem}
Suppose that a finite group $G$ acts on a smooth variety $V$ and
that $V/G$ has log-terminal singularities. Does $\M_{\st}(V/G)$
have the rationality like \eqref{eq:rat ex}? If $V$ is proper, then
does the quotient variety satisfy the Poincaré duality \eqref{eq:Poincare}?
\end{problem}

We can ask the same question, more generally, for varieties
with log-terminal singularities. However, if we look for pathological
phenomena, especially, a counterexample to resolution of singularities,
then quotient varieties associated with non-linear actions might be
good places to search (cf.~\cite{kerz2013cohomological}). 

\section{Generalization in various directions\label{sec:Generalization-in-various}}
\begin{problem}
\label{prob:motivic-mixed}Develop motivic integration for Deligne-Mumford
stacks over a complete discrete valuation ring of mixed characteristics
and apply it to the proof of McKay correspondence over such a
ring. 
\end{problem}

Sebag \cite{sebag2004integration} (see also \cite{chambert-loir2018motivic})
developed motivic integration for (formal) schemes over a complete
discrete valuation ring, which may have mixed characteristic like
$p$-adic integers $\ZZ_{p}$. In \cite{yasuda2004twisted,yasuda2006motivic,yasuda2020motivic},
the author developed motivic integration for (formal) Deligne-Mumford
stacks over the power series ring $k\tbrats$. The wild McKay correspondence
also holds over $k\tbrats$. It is natural to look for a further
generalization to the complete discrete valuation ring possibly of mixed
characteristic. Note that the $p$-adic integration in a similar setting
was developed in \cite{yasuda2017thewild} (cf.~\cite{groechenig2020mirrorsymmetry}),
which includes the case of mixed characteristic. The remaining nontrivial
part appears to be the construction of moduli spaces used in the theory,
in particular, the moduli space $\Delta_{G}$ of $G$-torsors over
$\Spec k\tpars$ or its variants (Problem \ref{prob:moduli mix}).
\begin{problem}
\label{prob:alg gp}Generalize the McKay correspondence to more general
groups/group schemes such as algebraic groups (of positive dimension),
abelian varieties, non-reduced finite group schemes, and pro-finite
groups. 
\end{problem}

The case of the group scheme $\alpha_{p}$ was discussed in \cite{tonini2019noteson}
(cf.~Problem \ref{prob:Delta general G}). This problem is closely
related to:
\begin{problem}
\label{prob:Artin}Generalize the motivic integration to Artin stacks. 
\end{problem}

Balwe \cite{balwe2015padicand} has worked in this direction,
but his results do not appear to be applicable to the context of the McKay correspondence.
A recent work by Satriano and Usatine \cite{satriano2020stringy}
also addresses the above problem with the view toward McKay-correspondence
type results. However, studies on Problems \ref{prob:alg gp} and \ref{prob:Artin}
are still at an early stage, and there is considerable work that remains.
\begin{problem}
Generalize the works \cite{groechenig2020mirrorsymmetry,loeser2021motivic}
on the Hitchin fibration and stringy invariants twisted by $\mu_{n}$-gerbes
to the wild case. 
\end{problem}

In \cite{groechenig2020mirrorsymmetry}, Groechenig-Wyss-Ziegler proved
the topological mirror symmetry conjecture of Hausel and Thaddeus
by means of $p$-adic integration on Deligne-Mumford stacks. This
conjecture is stated in terms of stringy Hodge numbers twisted by
$\mu_{n}$-gerbes. Later, this result was refined and proved at the
motivic level by Loeser-Wyss \cite{loeser2021motivic}. It would be
interesting to incorporate wild actions into research in this direction. 

\section{Relation to other theories on wild ramification\label{sec:Relation-to-other-wild}}

Wild ramification is an important subject in number theory and arithmetic
geometry. The reader is referred to \cite{saito2011wildramification,xiao2015ramification}
and references therein for more discussion on the subject. An approach
to wild ramification in dimensions $\ge2$ is to restrict ramification
to various curves. The wild McKay correspondence as well as the motivic
integration over wild Deligne-Mumford stacks can be viewed as yet
another theory concerning wild ramification in higher dimensions via
the restriction-to-curves approach. 
\begin{problem}
Relate the wild McKay correspondence and motivic integration over
wild Deligne-Mumford stacks with other theories on wild ramification. 
\end{problem}

For example, non-liftability of arcs along a wildly ramified finite
cover, which played important roles in works of arithmetic geometers
\cite{kato2016wildramification,kerz2009ondifferent}, also played
a crucial role in application of stringy invariants to the local étale
fundamental group \cite{carvajal-rojas2021onthe}. See Section
\ref{sec:quasi-etale}. 

The Grothendieck-Ogg-Shafarevich formula \cite[Exposé X]{1977cohomologie}
is a classical result on wild ramification, and there are many works
on generalizing it (see \cite{saito2011wildramification,xiao2015ramification}).
\begin{problem}
Can we show a Grothendieck-Ogg-Shafarevich type formula by using motivic
integration and/or a version of stringy invariant?
\end{problem}

\section{The derived wild McKay correspondence\label{sec:The-derived-wild}}
\begin{problem}
Generalize the McKay correspondence at the level of derived categories
as studied in \cite{kapranov2000kleinian,bridgeland2001themckay}
to the wild case. 
\end{problem}

Suppose that there exists a crepant resolution $Y\to\AA_{k}^{d}/G$.
The derived McKay correspondence denotes the equivalence 
\[
D^{b}(\mathrm{Coh}(Y))\cong D^{b}(\mathrm{Coh}^{G}(\AA_{k}^{d})),
\]
between the bounded derived category of coherent sheaves on $Y$ and the one of 
 equivariant coherent sheaves on $\AA_{k}^{d}$. This was
proved under some conditions in characteristic zero. The category
$\mathrm{Coh}^{G}(\AA_{k}^{d})$ always has infinite global dimension
in the wild case \cite{yi1994homological}, while it has global dimension
$d$ in the tame case. Therefore, the above equivalence never holds
in the wild case. We would need to find an alternative to $\mathrm{Coh}^{G}(\AA_{k}^{d})$
possibly from various notions of higher structures developed especially
in noncommutative algebraic geometry.

\section{The McKay correspondence over global fields\label{sec:Over-global-fields}}
\begin{problem}
Explore the McKay correspondence over global fields, in particular,
number fields. 
\end{problem}

The McKay correspondence in terms of stringy motive, equality \eqref{eq:wild},
can be regarded as the McKay correspondence over a local field, as
both sides of the equality have something to do with the ``local
field'' $k\tpars$. More precisely, the left side can be seen as
the volume of the set of integral points (i.e., $k\tbrats$-points)
on the quotient variety $\AA_{k\tbrats}^{d}/G$, while the right side
can be seen as the volume of the set of $G$-torsors over the punctured
formal disk. This point of view is more apparent in the point-counting
version of \eqref{eq:wild} (see \cite{wood2015massformulas,yasuda2017thewild}).
It is natural to look for a counterpart over a global field, in
particular, a number field. A prime candidate for such a theory would
concern an interplay of the \emph{Batyrev-Manin conjecture \cite{franke1989rational,batyrev1990surle}}
and the \emph{Malle conjecture} \cite{malle2004onthe,malle2002onthe},
as was discussed in \cite{yasuda2014densities,yasuda2015maninsconjecture}.
The former discusses the number of rational points with bounded height
on a Fano-type variety over a number field, while the latter discusses
the number of Galois extensions with bounded discriminant. Both conjectures
predict these numbers increase like
\begin{equation}
c\cdot B^{a}\cdot(\log B)^{b}\label{eq:growth}
\end{equation}
for some constants $c>0,a>0,b\ge0$, as the bound $B$ tends to infinity.
There are close relations among these constants in both contexts (see
\cite[Sec.\ 4]{yasuda2014densities} and \cite[Prop. \ 8.5]{yasuda2015maninsconjecture}),
which are also considered as versions of the McKay correspondence.
For example, in the Malle conjecture for a transitive subgroup $G\subset S_{n}$,
the exponent $a$ of $B$ is the reciprocal of the \emph{index}, denoted
by $\ind(G)$, induced from the $G$-action on $\{1,\dots,n\}$. The
index is related to the \emph{discrepancy} of a quotient variety,
a basic invariant of singularities \cite[Prop.\ 4.4]{yasuda2014densities}.
The constants $b$ in the two conjectures are related to the number
of $K$\emph{-conjugacy classes with minimal index} and to the number
of \emph{crepant divisors} over a quotient variety (see \cite[Conj.\ 5.6]{yasuda2015maninsconjecture})
respectively. They are related by \cite[Prop.\ 4.5]{yasuda2014densities}
(cf. \cite[Prop.\ 8.5]{yasuda2015maninsconjecture}). As for the constant
$c$, on the Malle side, we have the conjectural formula of Bhargava
for $G=S_{n}$ \cite[Conj.\ A]{bhargava2007massformulae}, which he
derived from his own formula for the local counterpart. On the Batyrev-Manin
side, we have Peyre's constant \cite{peyre1995hauteurs} (see also
\cite{peyre2003pointsde,batyrev1998tamagawa}), which is related to
a volume of adelic points. 
\begin{problem}
Prove one of the above two conjectures for some special case from
the known case of the other conjecture. 
\end{problem}

Le Rudulier's thesis \cite{lerudulier2014pointsalgebriques} is the
first such successful attempt. It would also be possible to interpret
works of Schmidt, Gao, and Guignard \cite{schmidt1995northcotts,gao1996onnorthcotts,guignard2017counting}
on counting algebraic points on a projective space as counting rational
points on a quotient variety via the correspondence of points (see
\cite[Sec.\ 4]{yasuda2014densities}). 
\begin{problem}
Count Galois extensions of a number field using a new weight function
derived from tuning modules. 
\end{problem}

Another importance of the McKay-correspondence viewpoint is to provide
new weights of the geometric origin for counting global fields, which
correspond to the $v$-function for local fields (see Section \ref{sec:The-v-function:-Fr=0000F6hlich=002019s}).
Dummit \cite{dummit2014counting} defined $\rho$-\emph{discriminant} as such a weight using the tuning module and studied the counting number
fields with bounded $\rho$-discriminant. Similarly, Yasuda \cite{yasuda2015maninsconjecture}
considered $V$-\emph{discriminant} and \emph{extended $V$-discriminant}
to discuss the relation between the Batyrev-Manin conjecture and the Malle
conjecture. Later, Ellenberg, Satriano, and Zureick-Brown \cite{ellenberg2021heights}
announced a conjecture on asymptotics of the number of rational points
on an algebraic stack with bounded height, which unifies the above
two conjectures. For this purpose, they introduced height functions
associated with vector bundles. In this new notion of height, the $v$-function
appears as the contribution of each non-Archimedean place (see also
\cite{landesman2021athesis}). Darda \cite{darda2021rational} also
defined height on stacks in connection with the conjectures of Batyrev-Manin
and Malle. The work \cite{odesky2020themoduli} of O'Desky and Rosen
is also considered as one aiming toward the same direction; they constructed
the moduli space of pairs of a Galois $\QQ$-algebra and a normal
element of it as an open subscheme of a certain quotient variety. 

\section{Stringy motives and quasi-étale Galois coverings\label{sec:quasi-etale}}

Let $Y\to X$ be a Galois quasi-étale (i.e., étale in codimension
1) cover of normal $\QQ$-Gorenstein varieties with only log-terminal
singularities. Suppose that the Galois group $G$ fixes some point
$y$ of $Y$. Then, we can define the $G$-quotient $\M_{\st}(Y)_{y}/G$
of $\M_{\st}(Y)_{y}$. If $x\in X$ denotes the image of $y$, then
we have the inequality
\begin{equation}
\M_{\st}(X)_{x}\ge\M_{\st}(Y)_{y}/G.\label{eq:q-et}
\end{equation}
Here, we compare these invariants by taking the Poincaré polynomial
realization and comparing their coefficients lexicographically. 
\begin{problem}
Is this inequality always strict unless the cover is trivial?
\end{problem}

In \cite{carvajal-rojas2021onthe}, the authors proved that the answer
is affirmative in dimension two as well as in dimension three with
the condition that either $Y$ has rational singularities or $G$
is a $p$-group. The difference between the two sides of \eqref{eq:q-et}
arises from the arcs on $X$ that do not lift to $Y$. To show the
strict inequality, we have to show that these non-liftable arcs form
a subset of the arc space of nonzero measure.
\begin{problem}
Consider a tower 
\[
X\leftarrow Y_{1}\leftarrow Y_{2}\leftarrow\cdots
\]
of Galois quasi-étale coverings as above and let $y_{i}\in Y_{i}$
and $G_{i}$ be the point and Galois group as above for the Galois
covering $Y_{i}\to X$. Does the sequence 
\[
\M_{\st}(Y_{1})_{y_{1}}/G_{1}\ge\M_{\st}(Y_{2})_{y_{1}}/G_{2}\ge\cdots
\]
satisfy the descending chain condition (DCC)?
\end{problem}

In the paper mentioned above, it is proved that the answer to this problem
is also affirmative in dimension two. By combining it with the strictness
of inequality, it can be observed that the local étale fundamental group of a
log-terminal surface singularity is finite (in any characteristic). 


\end{document}

%% file: macros.tex
\global\long\def\bigmid{\mathrel{}\middle|\mathrel{}}%

\global\long\def\AA{\mathbb{A}}%

\global\long\def\CC{\mathbb{C}}%

\global\long\def\FF{\mathbb{F}}%

\global\long\def\GG{\mathbb{G}}%

\global\long\def\LL{\mathbb{L}}%

\global\long\def\MM{\mathbb{M}}%

\global\long\def\NN{\mathbb{N}}%

\global\long\def\PP{\mathbb{P}}%

\global\long\def\QQ{\mathbb{Q}}%

\global\long\def\RR{\mathbb{R}}%

\global\long\def\SS{\mathbb{S}}%

\global\long\def\ZZ{\mathbb{Z}}%

\global\long\def\bA{\mathbf{A}}%

\global\long\def\ba{\mathbf{a}}%

\global\long\def\bb{\mathbf{b}}%

\global\long\def\bd{\mathbf{d}}%

\global\long\def\bf{\mathbf{f}}%

\global\long\def\bg{\mathbf{g}}%

\global\long\def\bh{\mathbf{h}}%

\global\long\def\bj{\mathbf{j}}%

\global\long\def\bm{\mathbf{m}}%

\global\long\def\bp{\mathbf{p}}%

\global\long\def\bq{\mathbf{q}}%

\global\long\def\br{\mathbf{r}}%

\global\long\def\bs{\mathbf{s}}%

\global\long\def\bt{\mathbf{t}}%

\global\long\def\bv{\mathbf{v}}%

\global\long\def\bw{\mathbf{w}}%

\global\long\def\bx{\boldsymbol{x}}%

\global\long\def\by{\boldsymbol{y}}%

\global\long\def\bz{\mathbf{z}}%

\global\long\def\bA{\mathbf{A}}%

\global\long\def\bB{\mathbf{B}}%

\global\long\def\bC{\mathbf{C}}%

\global\long\def\bD{\mathbf{D}}%

\global\long\def\bE{\mathbf{E}}%

\global\long\def\bF{\mathbf{F}}%

\global\long\def\bG{\mathbf{G}}%

\global\long\def\bM{\mathbf{M}}%

\global\long\def\bP{\mathbf{P}}%

\global\long\def\bS{\mathbf{S}}%

\global\long\def\bU{\mathbf{U}}%

\global\long\def\bV{\mathbf{V}}%

\global\long\def\bW{\mathbf{W}}%

\global\long\def\bX{\mathbf{X}}%

\global\long\def\bY{\mathbf{Y}}%

\global\long\def\bZ{\mathbf{Z}}%

\global\long\def\cA{\mathcal{A}}%

\global\long\def\cB{\mathcal{B}}%

\global\long\def\cC{\mathcal{C}}%

\global\long\def\cD{\mathcal{D}}%

\global\long\def\cE{\mathcal{E}}%

\global\long\def\cF{\mathcal{F}}%

\global\long\def\cG{\mathcal{G}}%

\global\long\def\cH{\mathcal{H}}%

\global\long\def\cI{\mathcal{I}}%

\global\long\def\cJ{\mathcal{J}}%

\global\long\def\cK{\mathcal{K}}%

\global\long\def\cL{\mathcal{L}}%

\global\long\def\cM{\mathcal{M}}%

\global\long\def\cN{\mathcal{N}}%

\global\long\def\cO{\mathcal{O}}%

\global\long\def\cP{\mathcal{P}}%

\global\long\def\cQ{\mathcal{Q}}%

\global\long\def\cR{\mathcal{R}}%

\global\long\def\cS{\mathcal{S}}%

\global\long\def\cT{\mathcal{T}}%

\global\long\def\cU{\mathcal{U}}%

\global\long\def\cV{\mathcal{V}}%

\global\long\def\cW{\mathcal{W}}%

\global\long\def\cX{\mathcal{X}}%

\global\long\def\cY{\mathcal{Y}}%

\global\long\def\cZ{\mathcal{Z}}%

\global\long\def\fa{\mathfrak{a}}%

\global\long\def\fb{\mathfrak{b}}%

\global\long\def\fc{\mathfrak{c}}%

\global\long\def\ff{\mathfrak{f}}%

\global\long\def\fj{\mathfrak{j}}%

\global\long\def\fm{\mathfrak{m}}%

\global\long\def\fp{\mathfrak{p}}%

\global\long\def\fs{\mathfrak{s}}%

\global\long\def\ft{\mathfrak{t}}%

\global\long\def\fx{\mathfrak{x}}%

\global\long\def\fv{\mathfrak{v}}%

\global\long\def\fD{\mathfrak{D}}%

\global\long\def\fJ{\mathfrak{J}}%

\global\long\def\fG{\mathfrak{G}}%

\global\long\def\fM{\mathfrak{M}}%

\global\long\def\fO{\mathfrak{O}}%

\global\long\def\fS{\mathfrak{S}}%

\global\long\def\fV{\mathfrak{V}}%

\global\long\def\fX{\mathfrak{X}}%

\global\long\def\fY{\mathfrak{Y}}%

\global\long\def\ru{\mathrm{u}}%

\global\long\def\rv{\mathbf{\mathrm{v}}}%

\global\long\def\rw{\mathrm{w}}%

\global\long\def\rx{\mathrm{x}}%

\global\long\def\ry{\mathrm{y}}%

\global\long\def\rz{\mathrm{z}}%

\global\long\def\a{\mathrm{a}}%

\global\long\def\AdGp{\mathrm{AdGp}}%

\global\long\def\Aff{\mathbf{Aff}}%

\global\long\def\Alg{\mathbf{Alg}}%

\global\long\def\age{\operatorname{age}}%

\global\long\def\Ann{\mathrm{Ann}}%

\global\long\def\Aut{\operatorname{Aut}}%

\global\long\def\B{\operatorname{\mathrm{B}}}%

\global\long\def\Bl{\mathrm{Bl}}%

\global\long\def\C{\operatorname{\mathrm{C}}}%

\global\long\def\calm{\mathrm{calm}}%

\global\long\def\center{\mathrm{center}}%

\global\long\def\characteristic{\operatorname{char}}%

\global\long\def\codim{\operatorname{codim}}%

\global\long\def\Coker{\mathrm{Coker}}%

\global\long\def\Conj{\operatorname{Conj}}%

\global\long\def\D{\mathrm{D}}%

\global\long\def\Df{\mathrm{Df}}%

\global\long\def\diag{\mathrm{diag}}%

\global\long\def\det{\operatorname{det}}%

\global\long\def\discrep#1{\mathrm{discrep}\left(#1\right)}%

\global\long\def\doubleslash{\sslash}%

\global\long\def\E{\operatorname{E}}%

\global\long\def\Emb{\operatorname{Emb}}%

\global\long\def\et{\textrm{ét}}%

\global\long\def\etop{\mathrm{e}_{\mathrm{top}}}%

\global\long\def\el{\mathrm{e}_{l}}%

\global\long\def\Exc{\mathrm{Exc}}%

\global\long\def\Fitt{\operatorname{Fitt}}%

\global\long\def\Gal{\operatorname{Gal}}%

\global\long\def\GalGps{\mathrm{GalGps}}%

\global\long\def\GL{\mathrm{GL}}%

\global\long\def\Grass{\mathrm{Grass}}%

\global\long\def\H{\operatorname{\mathrm{H}}}%

\global\long\def\hattimes{\hat{\times}}%

\global\long\def\hatotimes{\hat{\otimes}}%

\global\long\def\Hilb{\mathrm{Hilb}}%

\global\long\def\Hodge{\mathrm{Hodge}}%

\global\long\def\Hom{\operatorname{Hom}}%

\global\long\def\hyphen{\textrm{-}}%

\global\long\def\I{\operatorname{\mathrm{I}}}%

\global\long\def\id{\mathrm{id}}%

\global\long\def\Image{\operatorname{\mathrm{Im}}}%

\global\long\def\injlim{\varinjlim}%

\global\long\def\iper{\mathrm{iper}}%

\global\long\def\Iso{\operatorname{Iso}}%

\global\long\def\isoto{\xrightarrow{\sim}}%

\global\long\def\J{\operatorname{\mathrm{J}}}%

\global\long\def\Jac{\mathrm{Jac}}%

\global\long\def\Ker{\operatorname{Ker}}%

\global\long\def\Kzero{\operatorname{K_{0}}}%

\global\long\def\lcr{\mathrm{lcr}}%

\global\long\def\lcm{\operatorname{\mathrm{lcm}}}%

\global\long\def\length{\operatorname{\mathrm{length}}}%

\global\long\def\M{\operatorname{\mathrm{M}}}%

\global\long\def\MHS{\mathbf{MHS}}%

\global\long\def\mld{\mathrm{mld}}%

\global\long\def\mod#1{\pmod{#1}}%

\global\long\def\mRep{\mathbf{mRep}}%

\global\long\def\mult{\mathrm{mult}}%

\global\long\def\N{\operatorname{\mathrm{N}}}%

\global\long\def\Nef{\mathrm{Nef}}%

\global\long\def\nor{\mathrm{nor}}%

\global\long\def\NS{\mathrm{NS}}%

\global\long\def\op{\mathrm{op}}%

\global\long\def\ord{\operatorname{ord}}%

\global\long\def\P{\operatorname{P}}%

\global\long\def\PEff{\mathrm{PEff}}%

\global\long\def\PGL{\mathrm{PGL}}%

\global\long\def\pt{\mathbf{pt}}%

\global\long\def\pur{\mathrm{pur}}%

\global\long\def\perf{\mathrm{perf}}%

\global\long\def\pr{\mathrm{pr}}%

\global\long\def\Proj{\operatorname{Proj}}%

\global\long\def\projlim{\varprojlim}%

\global\long\def\Qbar{\overline{\QQ}}%

\global\long\def\R{\operatorname{\mathrm{R}}}%

\global\long\def\Ram{\operatorname{\mathrm{Ram}}}%

\global\long\def\rank{\operatorname{\mathrm{rank}}}%

\global\long\def\rig{\mathrm{rig}}%

\global\long\def\red{\mathrm{red}}%

\global\long\def\reg{\mathrm{reg}}%

\global\long\def\rep{\mathrm{rep}}%

\global\long\def\Rep{\mathbf{Rep}}%

\global\long\def\sbrats{\llbracket s\rrbracket}%

\global\long\def\Sch{\mathbf{Sch}}%

\global\long\def\sep{\mathrm{sep}}%

\global\long\def\Set{\mathbf{Set}}%

\global\long\def\sing{\mathrm{sing}}%

\global\long\def\sm{\mathrm{sm}}%

\global\long\def\SL{\mathrm{SL}}%

\global\long\def\Sp{\operatorname{Sp}}%

\global\long\def\Spec{\operatorname{Spec}}%

\global\long\def\Spf{\operatorname{Spf}}%

\global\long\def\ss{\mathrm{ss}}%

\global\long\def\st{\mathrm{st}}%

\global\long\def\Stab{\operatorname{Stab}}%

\global\long\def\Supp{\operatorname{Supp}}%

\global\long\def\spars{\llparenthesis s\rrparenthesis}%

\global\long\def\Sym{\mathrm{Sym}}%

\global\long\def\tame{\mathrm{tame}}%

\global\long\def\tbrats{\llbracket t\rrbracket}%

\global\long\def\top{\mathrm{top}}%

\global\long\def\tors{\mathrm{tors}}%

\global\long\def\tpars{\llparenthesis t\rrparenthesis}%

\global\long\def\Tr{\mathrm{Tr}}%

\global\long\def\ulAut{\operatorname{\underline{Aut}}}%

\global\long\def\ulHom{\operatorname{\underline{Hom}}}%

\global\long\def\ulIso{\operatorname{\underline{{Iso}}}}%

\global\long\def\ulSpec{\operatorname{\underline{{Spec}}}}%

\global\long\def\Utg{\operatorname{Utg}}%

\global\long\def\Unt{\operatorname{Unt}}%

\global\long\def\Var{\mathbf{Var}}%